\documentclass{article}

\usepackage{color}
\definecolor{darkblue}{rgb}{0.0, 0.0, 0.55}

\PassOptionsToPackage{square,numbers,sort&compress}{natbib}
\usepackage[colorlinks=true, bookmarksopen,
            linkcolor={darkblue},
            anchorcolor={black},
            citecolor={darkblue},
            filecolor={magenta},
            menucolor={black},
            plainpages=false,pdfpagelabels,
            urlcolor={darkblue}]{hyperref}

\usepackage[preprint]{neurips_2020}
\usepackage{definitions}

\usepackage[utf8]{inputenc} 
\usepackage[T1]{fontenc}    
\usepackage{url}            
\usepackage{booktabs}       
\usepackage{amsfonts}       
\usepackage{nicefrac}       
\usepackage{microtype}      
\usepackage{amsmath}
\usepackage{bm}
\usepackage{graphicx}
\usepackage{chngcntr}
\counterwithin{table}{section}

\usepackage[vlined, ruled, linesnumbered]{algorithm2e}
\SetKwInOut{Input}{Input}
\SetKwInOut{Output}{Output}
\SetKwInOut{Params}{Data}
\SetKwProg{Function}{def}{:}{}
\DontPrintSemicolon

\title{An adaptive stochastic gradient-free approach for high-dimensional blackbox optimization}

\author{%
    Anton Dereventsov \\
    Lirio AI Research\\
    Lirio, LLC\\
    Knoxville, TN 37923 \\
    \texttt{adereventsov@lirio.com}
    \And
    Clayton G.~Webster\thanks{Behavioral Reinforcement and Learning Lab (BReLL), Lirio, LLC.}\\
    Lirio AI Research\\
    Lirio, LLC\\
    Knoxville, TN 37923 \\
    \texttt{cwebster@lirio.com}
    \And
    Joseph Daws,~Jr.\\
    Lirio AI Research\\
    Lirio, LLC\\
    Knoxville, TN 37923\\
    \texttt{jdaws@lirio.com}
}

\begin{document}

\maketitle

\begin{abstract}
In this work, we propose a novel adaptive stochastic gradient-free (ASGF) approach for solving high-dimensional nonconvex optimization problems based on function evaluations.
We employ a directional Gaussian smoothing of the target function that generates a surrogate of the gradient and assists in avoiding bad local optima by utilizing nonlocal information of the loss landscape.
Applying a deterministic quadrature scheme results in a massively scalable technique that is sample-efficient and achieves spectral accuracy.
At each step we randomly generate the search directions while primarily following the surrogate of the smoothed gradient.
This enables exploitation of the gradient direction while maintaining sufficient space exploration, and accelerates convergence towards the global extrema.
In addition, we make use of a local approximation of the Lipschitz constant in order to adaptively adjust the values of all hyperparameters, thus removing the careful fine-tuning of current algorithms that is often necessary to be successful when applied to a large class of learning tasks.
As such, the ASGF strategy offers significant improvements when solving high-dimensional nonconvex optimization problems when compared to other gradient-free methods (including the so called "evolutionary strategies'') as well as iterative approaches that rely on the gradient information of the objective function.
We illustrate the improved performance of this method by providing several comparative numerical studies on benchmark global optimization problems and reinforcement learning tasks.
\end{abstract}

\section{Introduction}\label{sec:intro}
In this work, we introduce a novel adaptive stochastic gradient-free (ASGF) approach for blackbox optimization.
The ASGF method achieves improved performance when compared with existing gradient-free optimization (GFO) schemes (including those sometimes referred to as evolutionary strategies~\citep{salimans2017evolution,Hansen_CMA,Hansen_Ostermeier_CMA_01,liu2019trust,Hamalainen-PPO}), as well iterative approaches that rely on the gradient information of the objective function on several challenging benchmarks problems.
This new technique is designed to alleviate some of the more difficult challenges associated with successful deployment of machine learning models for complex tasks, namely: (i) high-dimensionality; (ii) nonconvexity and extraneous local extrema; and (iii) extreme sensitivity to hyperparameters.

Data sets, models and network architectures have increased to gargantuan size and complexity, see, e.g., \citep{brown2020language,ILSVRC15}.
In these settings the efficiency of backpropagation and automatic differentiation is diminished.
In addition, the use of iterative approaches that rely on the gradient information of the objective function for training can lead to very poor performance, due to extraneous local optima of the loss landscape~\citep{MorseStanley16}.
Such methods also perform poorly on machine learning tasks with non-differentiable objective functions, e.g., reinforcement learning tasks.
GFO methods, i.e., those that rely solely on function evaluations, are well suited to address these issues, and recently significant advances have been made towards conquering these challenges~\citep{salimans2017evolution,Hansen_CMA,hansen2006cma,liu2019trust,CRSTW18,zhang2020scalable}.
These methods are sometimes referred to as ``evolutionary strategies (ES),'' which are a class of algorithms inspired by natural evolution~\citep{Wierstra_NES_14}.
Recently, interest in ES methods has been reinvigorated and has become a popular approach in several machine learning problems such as training neural networks~\citep{Such-GA17,MorseStanley16,Cui_NeurIPS18} and reinforcement learning~\citep{zhang2020accelerating,Khadka-CERL,SigaudStulp13,salimans2017evolution}.
As such, these methods are particularly useful when solving optimization problems related to nonconvex and nonsmooth objective functions.

However, there are several ingredients necessary to execute such methods successfully on blackbox optimization problems, including: the computation of efficient search directions; size of the smoothing radius; the number and type of samples used to evaluate the function; and parameters related to the iterative update and candidate solutions.
More importantly, the lack of generalization cannot be overlooked, notwithstanding the fact that solutions to most blackbox optimization problems require significant hyper-parameter tuning, renders the results when computing with such methods to be highly problem dependent.
In contrast, the ASGF approach can be successfully applied to a wide class of optimization problems using a fixed set of hyper-parameters.
Although more sophisticated methods than brute-force searches have been explored for identifying good hyperparameters, e.g., \citep{journals/jmlr/BergstraB12}, an approach that diminishes this process is highly desirable.

The major contribution of this paper is the design and application of ASGF to challenging nonconvex optimization problems.
In particular, our approach is designed to be adaptive, gradient-free, massively parallelizable and scalable, and easily tuned with simple choices of hyperparameters.
Moreover, ASGF accelerates convergence towards the global extrema by exploiting an innovative iterative procedure.
At each step we primarily follow the direction of an approximate gradient of a smoothed function, while still maintaining the domain exploration by randomly generated search directions.
In addition, we explore the local geometry of the objective function and then adaptively adjust the parameter values at each iteration.

\subsection{Related Works}\label{sec:related}
GFO approaches include a large class of techniques to optimize an objective function based only on function values.
For example, see~\citep{RiosSahinidis13,Larson_et_al_19} for a general review on topic.
A particular class of GFO algorithms, referred to as evolutionary strategies~\citep{Wierstra_NES_14}, combines Gaussian smoothing~\citep{NesterovSpokoiny15,FKM05} and random search techniques~\citep{Rechenberg_ESbook_1973,Maheswaranathan_GuidedES}, which have been applied to a large class of learning tasks including, e.g., reinforcement learning~\citep{muller2018challenges,zhang2020accelerating,Khadka-CERL,SigaudStulp13,salimans2017evolution}, as well as training and optimizing neural networks~\citep{Such-GA17,MorseStanley16,Cui_NeurIPS18}.
Moreover, these strategies have been shown to be competitive when compared to iterative approaches that rely on the gradient information of the objective function~\citep{Peters-Schaal08,Mania18,BCCS19,Schulman-PPO,pourchot2018cem}, and further improved by employing either adaptation~\citep{salimans2017evolution,Hansen_CMA,Hansen_Ostermeier_CMA_01,liu2019trust,Hamalainen-PPO} or orthogonal directional derivatives~\cite{CRSTW18,zhang2020scalable}.
Our work herein combines the advantages of both approaches in order to overcome several of the grand challenges associated with solving high-dimensional nonconvex optimization problems.

\section{Background}\label{sec:BG}
Our objective is to solve for the global extrema of a high-dimensional nonconvex objective function $f : \R^d \rightarrow \R$.
Without loss of generality we consider the unconstrained blackbox optimization problem, parameterized by a $d$-dimensional vector $x = (x_1, \ldots, x_d) \in \R^d$, i.e.,
\start\label{eq:opt}
    \min_{x \in \R^d} f(x). 
\finish
Throughout this effort we assume that $f(x)$ is only available by virtue of function evaluations, and the gradient $\nabla f(x)$ is inaccessible, thus~\eqref{eq:opt} is typically solved with a derivative- or gradient-free optimization (GFO) method~\citep{FGKM18,NesterovSpokoiny15,RiosSahinidis13,Larson_et_al_19,MWDS18}.

\subsection{Gaussian smoothing for GFO methods}\label{sec:GS4GFO}
Similarly to the so-called ``evolutionary strategies''~\citep{Hansen_Ostermeier_CMA_01,hansen2006cma,salimans2017evolution,CRSTW18}, we introduce the notion of Gaussian smoothing~\citep{Nesterov_book_2004,NesterovSpokoiny15} of the objective function $f(x)$ in~\eqref{eq:opt}.
Let $\sigma > 0$ be the smoothing parameter and denote by $f_\sigma(x)$ the Gaussian smoothing of $f$ with radius $\sigma$, i.e.,
\begin{equation*}
    f_\sigma(x) = \frac{1}{\pi^{\nicefrac{d}{2}}} \int_{\mathbb{R}^d} f(x + \sigma \epsilon) 
    \, e^{-\|\epsilon\|_2^2} \,\mathrm{d} \epsilon
    = \mathbb{E}_{\varepsilon \sim \mathcal{N}(0,\mathbb{I}_d)}
    \big[ f(x + \sigma\epsilon) \big].
\end{equation*}
We remark that $f_\sigma(x)$ preserves important features of the objective function including, e.g., convexity, the Lipschitz constant, and is always differential even when $f(x)$ is not.
Furthermore, since $\|f - f_\sigma\|$ can be bounded by the Lipschitz constant, problem~\eqref{eq:opt} can be replaced by the smoothed version, i.e., $\min_{x \in \R^d} f_\sigma(x)$ (see, e.g., \citep{zhang2020scalable} and references therein).
The gradient of $f_\sigma(x)$ can be computed as
\start\label{eq:f_s_grad}
	\nabla f_\sigma(x) 
	= \frac{2}{\sigma \pi^{\nicefrac{d}{2}}}
	    \int_{\mathbb{R}^d} \epsilon f(x + \sigma \epsilon) \, e^{-\|\epsilon\|_2^2} \, \mathrm{d}\epsilon
    = \frac{2}{\sigma} \mathbb{E}_{\epsilon \sim \mathcal{N}(0,\mathbb{I}_d)}
        \big[ \epsilon f(x + \sigma \epsilon) \big].
\finish
Then for $M \in \mathbb{N}_+$, $\left\{\epsilon_j\right\}_{j=1}^M \stackrel{\mathrm{iid}}{\sim} \mathcal{N}(0,1)$, and the learning rate $\lambda\in\R$, traditional GFO methods estimate \eqref{eq:f_s_grad} via Monte Carlo (MC) sampling and provide an iterative update to the state $x$, given by
\start\label{eq:MC_grad}
    \nabla f_\sigma(x) 
    \approx \frac{2}{\sigma M} \sum_{m=1}^M \epsilon_m f(x + \sigma \epsilon_m)
    \quad\text{and}\quad
    x_{i+1} = x_i - \frac{2\lambda}{\sigma M}\sum_{m=1}^M \epsilon_m f(x_i + \sigma\epsilon_m)
\finish
respectively~\citep{salimans2017evolution}.
The primary advantages of such GFO approaches is that they are easy to implement, embarrassingly parallelizable, and can be easily scaled to include a large number workers.
On the other hand, MC methods (see, e.g., \citep{Fishman_96}) suffer from slow convergence rates, proportional to $M^{-1/2}$, even though such rates are independent of dimension $d$.
Minor improvements could be expected with the use of quasi-MC sampling~\citep{Caflisch} or even sparse grid approximations~\citep{Nobile:2008wf,Nobile:2008uc}, however, the combination of high-dimensional domains and nonconvex objective functions makes all such GFO strategies only amenable to solving blackbox optimization problems in low to moderate dimensions.

{\bf Directional Gaussian smoothing via Monte Carlo sampling.} A promising attempt to improve the efficiency and accuracy of the MC gradient estimate~\eqref{eq:MC_grad} is to consider decoupling the problem~\eqref{eq:f_s_grad} along $d$ orthogonal directions~\citep{CRSTW18}.
The gradient can be estimated by virtue of, e.g., an antithetic orthogonal sampling, i.e., 
\s
    \nabla f_\sigma(x) 
    \approx \frac{1}{\sigma M} \sum_{j=1}^M \epsilon_j \big( f(x + \sigma \epsilon_j) - f(x - \sigma \epsilon_j) \big),
\f
where $\{\epsilon_j\}_{j=1}^M$ are marginally distributed as $\mathcal{N}(0,1)$, and the joint distribution of $\left\{\epsilon_j\right\}_{j=1}^M$ is defined as follows: if $M \le d$, then the vectors are conditioned to be orthogonal almost surely.
If $M > d$, then each consecutive set of $d$ vectors is conditioned to be orthogonal almost surely, with distinct sets of $d$ vectors remaining independent.
Using the orthogonal directions, as opposed to the MC directions (as in~\eqref{eq:MC_grad}) improves the overall performance when approximating~\eqref{eq:f_s_grad}, however, due to the MC approximation along each orthogonal direction hinders the convergence of such methods suffers as the dimension increases.

{\bf Directional Gaussian smoothing via Gauss-Hermite quadrature.} An efficient approach for computing the decoupled integrals with {\em spectral accuracy} is to employ one-dimensional Gauss-Hermite quadrature.
This can be accomplished by letting $\Xi := (\xi_1, \ldots, \xi_d)$ be an orthonormal basis in $\mathbb{R}^d$ and by computing  directional derivatives $\nicefrac{\partial f_\sigma(x)}{\partial \xi_j}$ of $f_\sigma$ at point $x$ in the direction $\xi_j$, estimated as
\[
	\frac{\partial f_\sigma(x)}{\partial \xi_j} \approx \nabla f_\sigma(x |\, \xi_j)
	:= \frac{2}{\sigma \sqrt{\pi}} \int_\mathbb{R} v f(x + \sigma v \xi_j) \, e^{-v^2} \,\mathrm{d}v.
\]
Then the directional derivatives $\nabla f_\sigma(x_i |\, \xi_j)$ can be computed via Gauss-Hermite quadrature with $m_j \ge 3$ quadrature points, i.e., 
\start\label{eq:df_gh_quad}
	\nabla f_\sigma(x |\, \xi_j) 
	\approx \frac{2}{\sigma \sqrt{\pi}} \sum_{m=1}^{m_j} w_m p_m f(x + \sigma p_m \xi_j),
\finish
where $p_m$ are the roots of the Hermite polynomial of degree $m_j$ and $w_m$ are the corresponding weights (see, e.g., \citep{Handbook}).
Once the directional derivatives are computed, the estimate of the gradient of the smoothed function $f_\sigma$ at point $x$ can be computed as
\start\label{eq:grad_sigma}
	\nabla f_\sigma(x) = \sum_{j=1}^d \nabla f_\sigma(x |\, \xi_j) \, \xi_j.
\finish
This approach is considered in~\citep{zhang2020scalable} for applications to nonconvex blackbox optimization and later in the context of RL tasks~\citep{zhang2020accelerating}.
However, as described in Section~\ref{sec:intro}, this technique requires significant hyper-parameter tuning, which necessitates the development of our fully adaptive stochastic gradient-free method strategy.

\section{The adaptive stochastic gradient-free method (ASGF)}\label{sec:ASGF_method}
Before going into the details, we roughly outline the general flow of the algorithm.
At the beginning of an iteration $i$ the search directions $\Xi$ and smoothing parameter $\sigma$ are used to compute derivatives $\nabla f_\sigma(x_i |\, \xi_j)$ along the directions $\xi_j$ according to~\eqref{eq:df_gh_quad} and the gradient surrogate $\nabla f_\sigma(x_i)$ is estimated by~\eqref{eq:grad_sigma}.
The learning rate $\lambda$ is then selected based on certain local properties of the objective function and the candidate minimizer is updated by a step of gradient descent: $x_{i+1} = x_i - \lambda \, \nabla f_\sigma(x_i)$.
Finally, we update the search directions $\Xi$ and smoothing parameter $\sigma$, and proceed to the next iteration.
In the following sections we describe in detail each part of the process.

\subsection{The ASGF algorithm}\label{sec:ASGF_algo}
A central feature of the ASGF approach is the selection of the search directions $\Xi$.
Although the directional smoothing described in Section~\ref{sec:GS4GFO} holds for any set of orthonormal vectors $\Xi$, the choice of the updates to $\Xi$ has a significant impact on the realization of the optimization process.
For instance, taking steps mainly in the direction of the gradient $\nabla f(x)$ (assuming it exists) results in a form of (batch) gradient descent, while distributing updates across random directions is more in a style of stochastic gradient descent.
In ASGF the directions $\Xi$ are chosen in a way that balances efficiency and exploration as follows: for a particular iteration $i$ the algorithm the first direction $\xi_1$ is set to be the current estimate of the gradient of $f_\sigma(x_i)$, while the other directions $\xi_2, \ldots, \xi_d$ are chosen to complement $\xi_1$ to a random orthonormal basis in $\mathbb{R}^d$, i.e.,
\begin{equation}\label{eq:dir_update}
	\xi_1 = \frac{\nabla f_\sigma(x_i)}{\|\nabla f_\sigma(x_i)\|_2},
	\quad
	\xi_2, \ldots, \xi_d \in \mathbb{R}^d \text{ are such that } \Xi \text{ is an orthonormal basis}.
\end{equation}
Such an approach naturally combines the efficiency of exploiting the gradient direction $\xi_1$ (from now on referred to as `main' direction) while retaining the exploration ability provided by the stochastic directions $\xi_2, \ldots, \xi_d$ (called `auxiliary' directions), which are generated randomly on each iteration.

By splitting search directions into a single `main' direction and a set of `auxiliary' directions, we can improve the computational efficiency of the approach by using a different number of quadrature points for each of the two classes of directions.
For the `main' direction, i.e., the gradient direction, $\xi_1$ we use an adaptive scheme for establishing a suitable number of quadrature points.
Specifically, we estimate $\nabla f_\sigma(x |\, \xi_1)$ via~\eqref{eq:df_gh_quad} with an increasing numbers of quadrature points $m \in \{3, 5, 7, \ldots\}$ until we obtain two estimates that differ less than some threshold $\varepsilon_m$.
Since the `auxiliary' directions $\xi_2, \ldots, \xi_d$ mainly serve an exploration role, a fixed small number of quadrature points is used.
In the numerical experiments presented in section~\ref{sec:numerics}, we use $\varepsilon_m = .1$ and $m_2 = \ldots = m_d = 5$ points.

Another key aspect of ASGF is the adaptive selection of the learning rate.
Instead of using a fixed value for the learning rate or a predetermined schedule, the geometry of the target function is used to derive the step size.
For each direction $\xi_j$ the values $\{f(x + \sigma p_m \xi_j)\}_{m=1}^{m_j}$, sampled in~\eqref{eq:df_gh_quad}, are used to estimate the directional local Lipschitz constants $L_j$ as
\begin{equation}\label{eq:lip_loc}
	L_j = \max_{1 \le m < m_j} 
	\left| \frac{f(x + \sigma p_{m+1} \xi_j) - f(x + \sigma p_m \xi_j)}{\sigma (p_{m+1} - p_m)} \right|.
\end{equation}
The learning rate is derived from the smoothing parameter $\sigma$ and a running average over Lipschitz constants along the `main' direction computed on previous iterations, denoted $L_\nabla$.
Namely,
\begin{equation}\label{eq:lr}
    L_\nabla \leftarrow (1 - \gamma_L) \, L_1 + \gamma_L \, L_\nabla
    \quad\text{and}\quad
    \lambda = \nicefrac{\sigma}{L_\nabla},
\end{equation}
where $0 \le \gamma_L < 1$ is a the smoothing factor (we use $\gamma_L = .9$).
The averaging of $L_\nabla$ is used to add ``momentum-like" dynamics to the learning rate schedule.
The construction of $\lambda$ ensures that the larger steps are taken in correspondence to stronger smoothing of the target function $f$.
Conversely, smaller steps are taken when the local geometry along the `main' direction is less regular.
We note that such approach for the learning rate selection is known in the optimization community, see, e.g., \citep{nesterov2018lectures}.
In practice, we observe that such a scheme is efficient in a wide range of optimization settings as shown in Section~\ref{sec:numerics}.

Lastly, another substantial advantage of ASGF is the massive scalability since the computation of the directional derivatives~\eqref{eq:df_gh_quad}, as well as the local Lipschitz constants~\eqref{eq:lip_loc}, may be distributed in parallel across as many as $m_1 + \ldots + m_d$ workers.
Since the only communication required between parallel workers is the value of the target function at the quadrature points assigned to them, tremendous speed up is observed in practice when implemented with parallelism.

\begin{algorithm}
	\caption{ASGF algorithm}\label{alg:asgf}
	\Input{ function $f$, initial state $x_0$, initial smoothing parameter $\sigma_0$, termination tolerance $\varepsilon_x$}
	\Output{ point of minimum $\widetilde{x}$}
	\BlankLine
	set $\widetilde{x} \leftarrow x_0$ and $\sigma \leftarrow \sigma_0$, and let $\Xi$ be a random orthonormal basis\;
	\For{$i \leftarrow 0$ \KwTo \tt{maxiter}}{
		\For{$j \leftarrow 1$ \KwTo $d$ (in parallel)}{
			compute derivatives $\nabla f_\sigma(x_i |\, \xi_j)$ by~\eqref{eq:df_gh_quad} and 
			estimate local Lipschitz constants $L_j$ by~\eqref{eq:lip_loc}\;
		}
		average Lipschitz constant $L_\nabla$ and compute learning rate $\lambda$ by~\eqref{eq:lr}\;
		assemble the gradient $\nabla f_\sigma(x_i)$ via~\eqref{eq:grad_sigma} and 
		update the minimizer $x_{i+1} = x_i - \lambda \, \nabla f_\sigma(x_i)$\;
		\If{$f(x_{i+1}) < f(\widetilde{x})$}{
			update the best state $\widetilde{x} \leftarrow x_{i+1}$\;
		}
		\uIf{$\|x_{i+1} - x_i\|_2 < \varepsilon_x$}{
			\bf{break}\;
		}\Else{
			update smoothness parameter $\sigma$ and search directions $\Xi$ by Algorithm~\ref{alg:par_update}\;
		}
	}
\end{algorithm}

\subsection{Adaptivity of the smoothing parameter}
The updates in the smoothing parameter are made with respect to the local geometry of the function $f$ which is characterized be the constants $L_1, \ldots, L_d$, estimated by~\eqref{eq:lip_loc}.
To update $\sigma$, the values of directional derivatives $\nabla f_\sigma(x |\, \xi_j)$ are compared to the corresponding local Lipschitz constant $L_j$.
If the ratio $|\nicefrac{\nabla f_\sigma(x |\, \xi_j)}{L_j}|$ is sufficiently small, the value of $\sigma$ is decreased, and vice versa, when the ratio is large, $\sigma$ is increased.

To address the tendency of methods for nonconvex minimization to become stuck in a local minimum, in ASGF we `reset parameters' by `forgetting' the gathered information about the loss landscape under certain conditions.
Whenever smoothing parameter becomes sufficiently small ($\sigma < \rho \, \sigma_0$, where $0 <\rho < 1$ is a specified threshold), the values of $\sigma$ and $\Xi$ are reset to their initial values, which typically allows the optimizer to escape the particular neighborhood it is currently in.
Additionally, since frequent resets would be undesirable, we control the maximum number of allowed resets by a parameter $r$ (we use $r = 2$ and $\rho = .01$).
We note that similar approaches to escaping local minima are well-known in machine learning community, e.g.~\citep{auger2005restart} and~\cite{eriksson2019scalable}.

\begin{algorithm}
	\caption{Parameter update}\label{alg:par_update}
	\Input{ smoothing parameter $\sigma$, gradient $\nabla f_\sigma(x_i)$, local Lipschitz constants $L_1, \ldots, L_d$}
	\Params{ number of resets $r$ and reset factor $\rho$, decay rate $\gamma_\sigma$, threshold parameters $A, B$ and their change rates $A_+, A_-, B_+, B_-$}
	\Output{ smoothness parameter $\sigma$ and directions $\Xi$}
	\BlankLine
	\uIf{$r > 0$ {\bf and} $\sigma < \rho \, \sigma_0$}{
		assign $\Xi$ to be a random orthonormal basis and set $\sigma \leftarrow \sigma_0$\;
		set $A, B$ to their initial values and
		change number of resets $r \leftarrow r - 1$\;
	}\Else{
		update search directions $\Xi$ by~\eqref{eq:dir_update}\;
		\uIf{$\max_{1 \le j \le d} |\nicefrac{\nabla f_\sigma(x |\, \xi_j)}{L_j}| < A$}{
			decrease smoothing $\sigma \leftarrow \sigma * \gamma_\sigma$ and 
			lower threshold $A \leftarrow A * A_-$\;
		}\uElseIf{$\max_{1 \le j \le d} |\nicefrac{\nabla f_\sigma(x |\, \xi_j)}{L_j}| > B$}{
			increase smoothing $\sigma \leftarrow \sigma / \gamma_\sigma$ and 
			upper threshold $B \leftarrow B * B_+$\;
		}\Else{
			increase lower threshold $A \leftarrow A * A_-$ and
			decrease upper threshold $B \leftarrow B * B_+$\;
		}
	}
\end{algorithm}

\section{Numerical experiments}\label{sec:numerics}
In this section we demonstrate the empirical performance of ASGF on a range of established blackbox optimization problems.
Specifically, we consider the optimization benchmarks from the Virtual Library of Simulation Experiments\footnote{\url{https://www.sfu.ca/~ssurjano/optimization.html}} with the standard choice of parameter values and domains.
In Section~\ref{sec:low_dim} we compare ASGF to other algorithms for nonconvex optimization in the low-dimensional setting and in Section~\ref{sec:high_dim} we showcase the performance of ASGF in the high-dimensional setting.

We note that, unlike many existing algorithms, ASGF does not require careful selection of hyperparameters in order to obtain state-of-the-art results since the parameters are adjusted throughout the realization of the algorithm.
As such, we use the exact same set of hyperparameters for all of the stated examples.
This showcases the adaptability and stability of our method to the hyperparameters selection, which is an essential feature in the blackbox setting.
We would also like to point out that with minor tweaking ASGF can achieve better results for each of the presented examples, however the main purpose of this section is to showcase the adaptive nature of the ASGF to automatically determine the suitable parameter values and that we abstain from any kind of hyperparameter tuning.
Namely, we use the following set of parameters for ASGF: $\gamma_s=.9$, $m=5$, $A=.1$, $B=.9$, $A_-=.95$, $A_+=1.02$, $B_-=.98$, $B_+=1.01$, $\gamma_L=.9$, $r=2$, $\rho=.01$, $\varepsilon_m=.1$, $\varepsilon_x=10^{-6}$.
For the initial value of $\sigma$ we use the heuristic $\sigma_0 = \nicefrac{\operatorname{diam}(\Omega)}{10}$, where $\Omega \subset \mathbb{R}^d$ is the spatial domain from which the initial state $x_0$ is sampled.
Nevertheless, due to the adaptive design of ASGF, the values of the above hyperparameters could be changed without significant distinction in the resulting performance.
All of the presented numerical experiments are performed in Python and the source code reproducing the stated results is publicly available at~\url{https://github.com/joedaws/ASGF}.

\subsection{Low-dimensional optimization}\label{sec:low_dim}
Even though ASGF is designed with a high-dimensional setting in mind, in order to provide an extensive comparison with other methods, in this section we consider a wide range of optimization benchmarks, presented in Table~\ref{tab:low_dim}.
Here we compare the following algorithms: ASGF (ours), Directional Gaussian Smoothing (DGS, see~\citep{zhang2020scalable}), and Covariance Matrix Adaptation (CMA, see~\citep{hansen2006cma}). 

\begin{table}[htbp]
	\label{tab:low_dim}
	\centering\small
	\caption{Convergence, average number of iterations and function evaluations over $100$ simulations.}
	\begin{tabular}{l||rrr|rrr|rrr}
		& \multicolumn{3}{c}{Convergence} & \multicolumn{3}{|c|}{Iterations} & \multicolumn{3}{c}{Function evaluations}
		\\\cline{2-10}
		Benchmark & ASGF & DGS & CMA & ASGF & DGS & CMA & ASGF & DGS & CMA
		\\\hline\hline
		Branin & 100\% & 64\% & 100\% & 335 & 2,568 & 70 & 3,820 & 23,113 & 418
		\\
		Cross-in-Tray & 99\% & 76\% & 79\% & 335 & 3,939 & 74 & 9,031 & 98,476 & 446
		\\
		Dropwave & 100\% & 99\% & 1\% & 524 & 1,232 & 77 & 44,645 & 40,646 & 463
		\\
		Sphere 10d & 100\% & 100\% & 100\% & 16 & 50 & 197 & 669 & 2,066 & 1,969
		\\
		Ackley 2d & 95\% & 97\% & 86\% & 78 & 1,202 & 77 & 3,774 & 10,816 & 462
		\\
		Ackley 5d & 100\% & 99\% & 91\% & 54 & 2,147 & 151 & 2,703 & 45,100 & 1,208
		\\
		Ackley 10d & 100\% & 94\% & 97\% & 56 & 3,068 & 233 & 3,582 & 125,751 & 2,333
		\\
		Levy 2d & 100\% & 100\% & 97\% & 352 & 1,704 & 72 & 5,043 & 56,231 & 430
		\\
		Levy 5d & 100\% & 100\% & 81\% & 475 & 760 & 146 & 12,909 & 61,583 & 1,172
		\\
		Levy 10d & 100\% & 100\% & 61\% & 456 & 714 & 240 & 22,353 & 114,998 & 2,404
		\\
		Rastrigin 2d & 96\% & 100\% & 14\% & 86 & 2,078 & 83 & 2,785 & 85,179 & 496
		\\
		Rastrigin 5d & 100\% & 100\% & 0\% & 2,075 & 1,788 & $-$ & 159,564 & 180,539 & $-$
		\\
		Rastrigin 10d & 100\% & 100\% & 0\% & 2,430 & 1,562 & $-$ & 232,258 & 314,035 & $-$
	\end{tabular}
\end{table}

The presented experiments are performed over $100$ independent simulations.
All the algorithms start a simulation with the same initial guess $x_0$ sampled at random from the spatial domain $\Omega \subset \mathbb{R}^d$.
A simulation is considered successful if an algorithm returns a minimizer that achieves a value within $10^{-4}$ of the global minimum.
For ASGF the hyperparameters are the same across all the examples and are stated in the preamble of Section~\ref{sec:numerics}.
For DGS we perform a hyperparameter search over the following grid: $\lambda \in \{.001, .003, .01, .03, .1\}$, $m \in \{5, 9, 13, 17, 21\}$, $\sigma \in \{1, 2, 3, 4, 5, 6, 7, 8, 9, 10\}$, $\gamma \in \{.001, .003, .01, .03, .1\}$ and report only the best obtained results for every example.
For CMA the value of $\sigma$ is set to be the same as for DGS in the same setting.
Lastly, the comparison of ASGF to additional algorithms (such as BFGS, Nelder--Mead, and Powell) is provided in Appendix~\ref{sec:appendix}, however the most competitive algorithms (namely, DGS and CMA) are presented in Table~\ref{tab:low_dim}.

\subsection{High-dimensional optimization}\label{sec:high_dim}
\begin{figure}[htbp]
	\includegraphics[width=.49\linewidth]{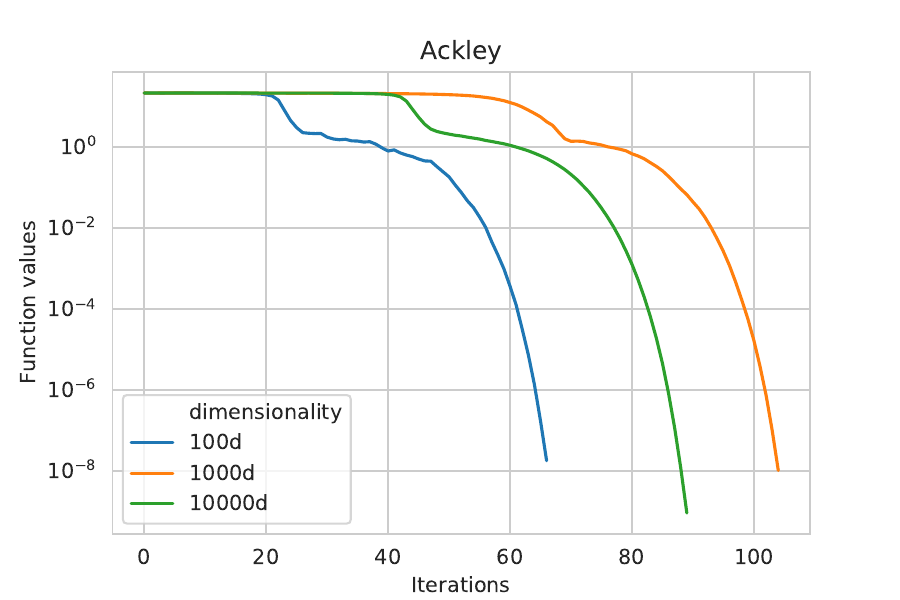}
	\includegraphics[width=.49\linewidth]{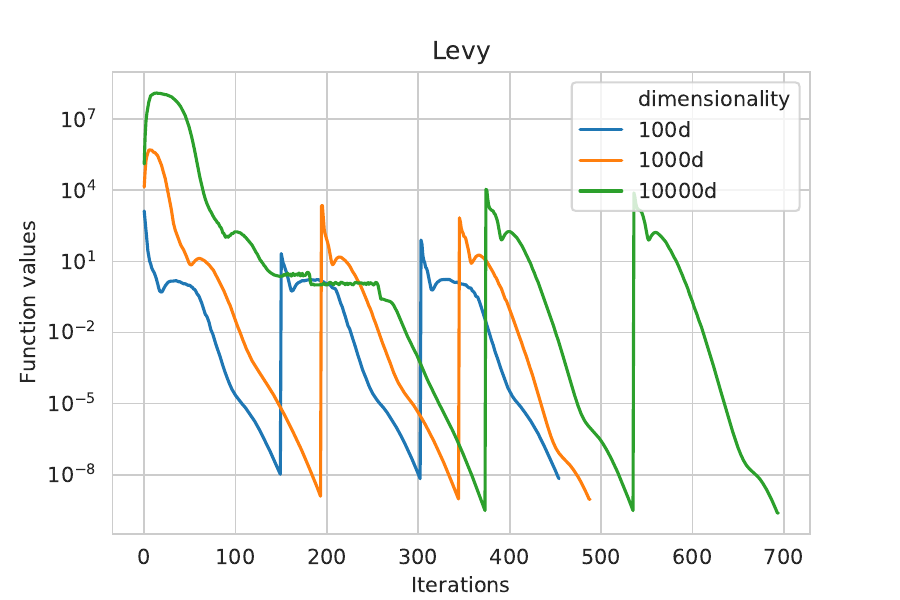}
	\\
	\includegraphics[width=.49\linewidth]{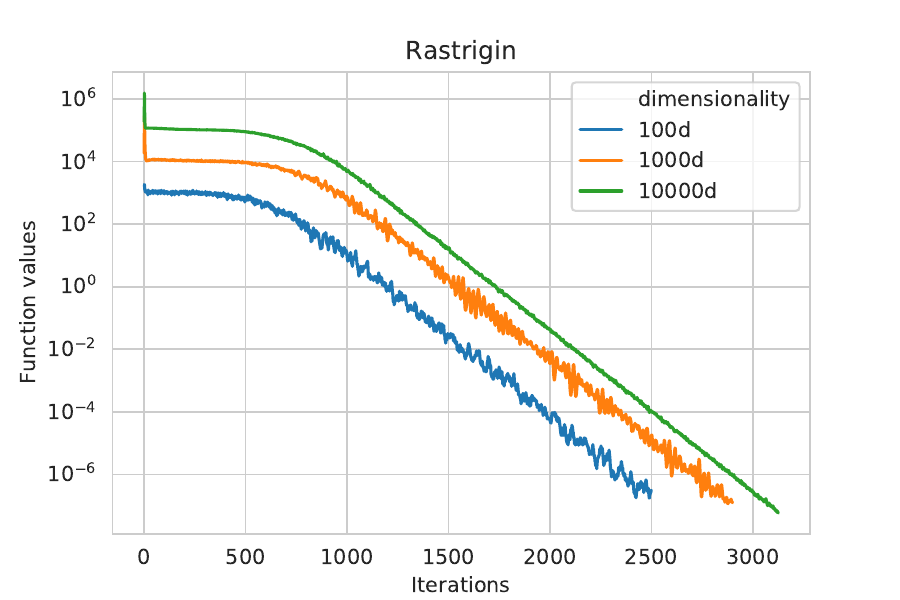}
	\includegraphics[width=.49\linewidth]{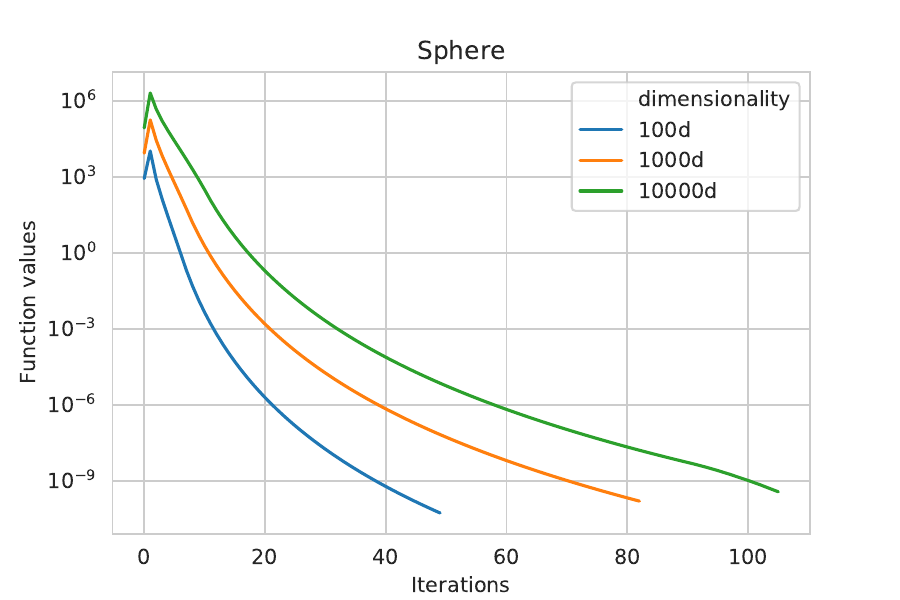}	
    \caption{Performance of ASGF on high-dimensional test functions.}
    \label{fig:high_dim}
\end{figure}

In this section we demonstrate the performance of ASGF in the setting of high-dimensional blackbox optimization on $100$, $1000$, and $10000$-dimensional benchmarks (namely, Ackley, Levy, Rastrigin, and Sphere functions).
All the simulations converged successfully, regardless of the initial state $x_0$, and the particular optimization trajectories are displayed in Figure~\ref{fig:high_dim}.
The average numbers of iterations and function evaluations are given in Table~\ref{tab:high_dim}.
The hyperparameters used for ASGF are the same ones as in Section~\ref{sec:low_dim} and are specified in the preamble of Section~\ref{sec:numerics}.

We note that the irregularity of the optimization trajectories in Figure~\ref{fig:high_dim} in case of Rastrigin function is caused by the highly oscillatory nature of the loss landscape.
The irregularity on Levy benchmark is due to the `parameter reset' feature of ASGF.
In particular, it can be observed that the `parameter reset' indeed helps the algorithm to escape local optima and converge to the global minimum.
In cases of Ackley and Sphere functions, due to distinguished geometry of the global minimum, the parameter reset is not triggered, hence the optimization trajectories are smooth.

\begin{table}[htbp]
    \label{tab:high_dim}
	\centering\small
	\caption{Average number of iterations and function evaluations of ASGF.}
	\begin{tabular}{l|r|r||l|r|r}
		Benchmark & Iterations & Evaluations & Benchmark & Iterations & Evaluations
		\\\hline\hline
        Ackley 100d & 66 & 27,343 & Rastrigin 100d & 2,995 & 1,290,215
        \\
        Ackley 1000d & 103 & 414,298 & Rastrigin 1000d & 2,901 & 11,625,963
        \\
        Ackley 10000d & 89 & 3,548,775 & Rastrigin 10000d & 3,206 & 114,845,172
        \\
        Levy 100d & 452 & 184,176 & Sphere 100d & 48 & 19,381
        \\
        Levy 1000d & 508 & 2,037,076 & Sphere 1000d & 76 & 303,508
        \\
        Levy 10000d & 617 & 24,698,140 & Sphere 10000d & 112 & 4,480,337
	\end{tabular}
\end{table}

\subsection{Reinforcement learning tasks}\label{sec:rl_examples}
The episodic reinforcement learning problem is formalized as a discrete control process over time steps where an agent interacts with its environment $\mathcal{E}$.
At each time step $t \in \{1,\ldots T\}$ the agent is presented with a state $s_t \in \mathcal{S}$ and correspondingly takes an action $a_t \in \mathcal{A}$.
This results in the agent receiving a reward $r_t = r(s_t,a_t) \in \mathcal{R}$ and the state transitions to $s_{t+1} \in \mathcal{S}$ by interfacing with the environment, i.e., $s_{t+1} = \mathcal{E}(s_t,a_t)$.
The agent achieves its goal by learning a policy $\pi: \mathcal{S} \rightarrow \mathcal{A}$ that maximizes the cumulative return objective function $f(\pi)$.
In machine learning, the agent uses a policy parameterized by a neural network, i.e., $\pi_x$, where the parameters 
$x=(x_1,\ldots, x_d)\in\R^d$ represent the weights of a neural network.  The architecture used in the example below is a two-layer network with 8 nodes per layer and the {\tt Tanh} activation function.
A good policy $\pi_{x}$ is obtained by solving the following nonconvex high-dimensional blackbox optimization problem 
\start\label{eq:RL}
    \max_{x \in \R^d} f(\pi_{x}),
    \quad\text{where}\quad
    f(\pi) := \sum_{t=1}^T \mathbb{E}_{\pi} \big[ r(s_t, a_t) \big],
\finish
where $a_t$ is sampled from the policy $\pi$.
In the numerical examples the initial state $s_0$ is chosen at random. 

\begin{figure}[htbp]
	\includegraphics[width=.49\linewidth]{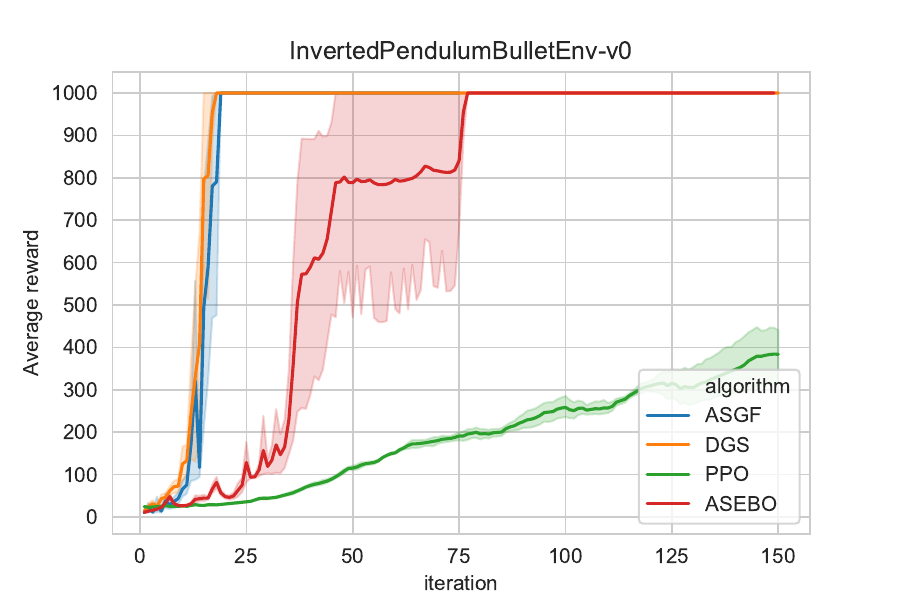}
	\includegraphics[width=.49\linewidth]{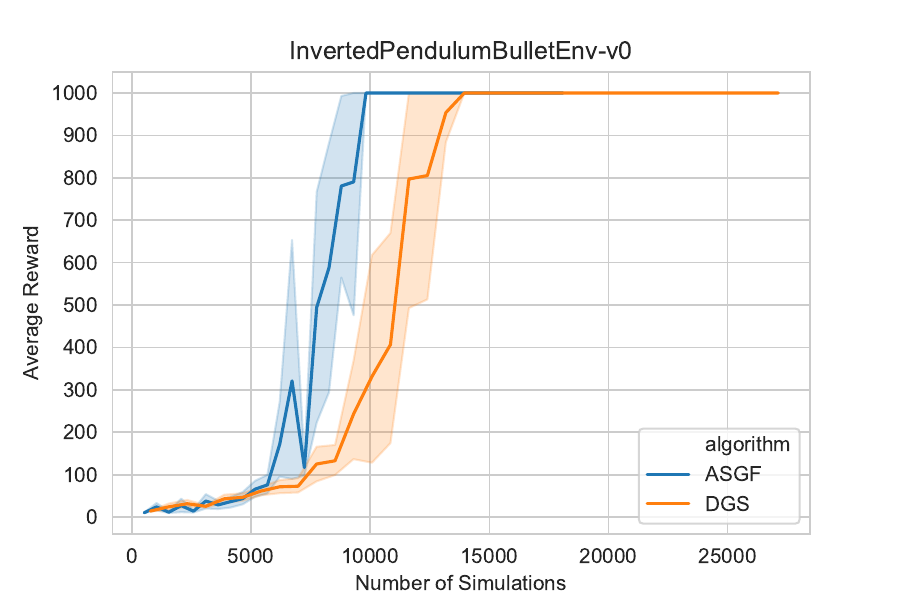}
    \caption{Comparison of several algorithms on {\tt InvertedPendulumBulletEnv-v0}.}
    \label{fig:rl}
\end{figure}
The performance of PPO~\citep{Hamalainen-PPO,Schulman-PPO}, ASEBO~\citep{Choromanski_ES-Active}, DGS~\cite{zhang2020scalable}, and ASGF is compared in Figure~\ref{fig:rl} for the {\tt InvertedPendulumBulletEnv-v0} reinforcement learning environment from the PyBullet library~\citep{coumans2016pybullet} using the OpenAI Gym library~\citep{brockman2016openai} interface where
the default episode length of $1000$ is used.
In each of the plots, the vertical axis is the average reward over $4$ training runs of each of the algorithms using different random seeds.
The horizontal axis of the figure on the left is labelled iterations where one iteration is one update of the network parameters.
It is reasonable to compare these methods using this metric since the computations between these updates is able to be parallelized.
Notice that ASGF and DGS both achieve good performance in much fewer iterations than PPO and ASEBO.
The horizontal axis of the plot on the right is the total number of simulations required to obtain a certain average reward.
Since ASGF uses an adaptive number of quadrature points, it needs fewer simulations to obtain a good average return.

\section{Conclusions}\label{sec:conclusions}
In this work we introduce an adaptive stochastic gradient-free method designed for solving high-dimensional nonconvex blackbox optimization problems.
The combination of hyperparameter adaptivity, massive scalability, and relative ease of implementation, makes ASGF prominent method for many practical applications.
The presented numerical examples empirically confirm that our method avoids many of the common pitfalls and overcomes challenges associated with high-dimensional nonconvex optimization.

Despite the successful demonstration of ASGF on several benchmark optimization problems in Section~\ref{sec:numerics}, we acknowledge that a single demonstration in the reinforcement learning domain is not a convincing case of the efficacy of ASGF for reinforcement learning, but we present it more as a proof of concept rather than the claim of superiority.
Our primary objective was to ensure that ASGF could successfully outperform existing GFO methods.
However, we feel strongly that since ASGF enables exploitation of the gradient direction while maintaining sufficient space exploration, we will also be able to accelerate the convergence on several more complicated RL tasks.
This is certainly the direction we are working on now and will be the focus of our future efforts.

\bibliography{references}

\newpage
\appendix
\section{Comparison of optimization algorithms in low-dimensional setting}\label{sec:appendix}
In this addendum we expand on the results of Section~\ref{sec:low_dim} and provide a more detailed comparison of different optimization methods, presented in Tables~A.1--A.3.
Specifically, we compare our Adaptive Stochastic Gradient-Free algorithm (ASGF), Directional Gaussian Smoothing (DGS, \citep{zhang2020scalable}), Covariance Matrix Adaptation (CMA, \citep{hansen2006cma}), Powell's conjugate direction method (Powell, \citep{powell1964efficient}), Nelder--Mead simplex direct search method (Nelder--Mead, \citep{nelder1965simplex}), and Broyden--Fletcher--Goldfarb--Shanno quasi-Newton method (BFGS, \citep{nocedal2006numerical}).
The hyperparameter choice for the relevant algorithms is discussed in Section~\ref{sec:low_dim}.
Each algorithm is tested on $100$ randomly sampled initial states, which are identical across all algorithms.
The simulation counts as successful if the returned minimizer achieves a function value that is within $10^{-4}$ of the global minimum.

\begin{table}[htbp]
	\centering\small
	\label{tab:low_dim_conv}
	\caption{Success rate of algorithms in terms of convergence to global minimum.}
	\begin{tabular}{l|rrrrrr}
		& ASGF & DGS & CMA & Powell & Nelder--Mead & BFGS
		\\\hline
		Branin & 100\% & 64\% & 100\% & 100\% & 100\% & 100\%
		\\
		Cross-in-Tray & 99\% & 76\% & 79\% & 14\% & 10\% & 10\%
		\\
		Dropwave & 100\% & 99\% & 1\% & 2\% & 0\% & 3\%
		\\
		Sphere 10d & 100\% & 100\% & 100\% & 100\% & 100\% & 100\%
		\\
		Ackley 2d & 95\% & 97\% & 86\% & 27\% & 1\% & 0\%
		\\
		Ackley 5d & 100\% & 99\% & 91\% & 2\% & 0\% & 0\%
		\\
		Ackley 10d & 100\% & 94\% & 97\% & 0\% & 0\% & 0\%
		\\
		Levy 2d & 100\% & 100\% & 97\% & 13\% & 7\% & 11\%
		\\
		Levy 5d & 100\% & 100\% & 81\% & 1\% & 0\% & 1\%
		\\
		Levy 10d & 100\% & 100\% & 61\% & 0\% & 0\% & 0\%
		\\
		Rastrigin 2d & 96\% & 100\% & 14\% & 33\% & 1\% & 3\%
		\\
		Rastrigin 5d & 100\% & 100\% & 0\% & 6\% & 0\% & 0\%
		\\
		Rastrigin 10d & 100\% & 100\% & 0\% & 1\% & 0\% & 0\%
	\end{tabular}
\end{table}

\begin{table}[htbp]
	\centering\small
	\label{tab:low_dim_itr}
	\caption{Average number of iterations on successful simulations.}
	\begin{tabular}{l|rrrrrr}
		& ASGF & DGS & CMA & Powell & Nelder--Mead & BFGS
		\\\hline
		Branin & 335 & 2,568 & 70 & 4 & 60 & 8
		\\
		Cross-in-Tray & 335 & 3,939 & 74 & 3 & 56 & 9
		\\
		Dropwave & 524 & 1,232 & 77 & 3 & $-$ & 5
		\\
		Sphere 10d & 16 & 50 & 197 & 2 & 1,972 & 4
		\\
		Ackley 2d & 78 & 1,202 & 77 & 6 & 79 & $-$
		\\
		Ackley 5d & 54 & 2,147 & 151 & 6 & $-$ & $-$
		\\
		Ackley 10d & 56 & 3,068 & 233 & $-$ & $-$ & $-$
		\\
		Levy 2d & 352 & 1,704 & 72 & 4 & 64 & 12
		\\
		Levy 5d & 475 & 760 & 146 & 10 & $-$ & 20
		\\
		Levy 10d & 456 & 714 & 240 & $-$ & $-$ & $-$
		\\
		Rastrigin 2d & 86 & 2,078 & 83 & 4 & 65 & 11
		\\
		Rastrigin 5d & 2,075 & 1,788 & $-$ & 5 & $-$ & $-$
		\\
		Rastrigin 10d & 2,430 & 1,562 & $-$ & 5 & $-$ & $-$
	\end{tabular}
\end{table}

\begin{table}[htbp]
	\centering\small
	\label{tab:low_dim_feval}
	\caption{Average number of function evaluations on successful simulations.}
	\begin{tabular}{l|rrrrrr}
		& ASGF & DGS & CMA & Powell & Nelder--Mead & BFGS
		\\\hline
		Branin & 3,820 & 23,113 & 418 & 126 & 116 & 38
		\\
		Cross-in-Tray & 9,031 & 98,476 & 446 & 114 & 109 & 40
		\\
		Dropwave & 44,645 & 40,646 & 463 & 117 & $-$ & 57
		\\
		Sphere 10d & 669 & 2,066 & 1,969 & 179 & 2,778 & 78
		\\
		Ackley 2d & 3,774 & 10,816 & 462 & 327 & 151 & $-$
		\\
		Ackley 5d & 2,703 & 45,100 & 1,208 & 700 & $-$ & $-$
		\\
		Ackley 10d & 3,582 & 125,751 & 2,333 & $-$ & $-$ & $-$
		\\
		Levy 2d & 5,043 & 56,231 & 430 & 111 & 123 & 69
		\\
		Levy 5d & 12,909 & 61,583 & 1,172 & 712 & $-$ & 196
		\\
		Levy 10d & 22,353 & 114,998 & 2,404 & $-$ & $-$ & $-$
		\\
		Rastrigin 2d & 2,785 & 85,179 & 496 & 150 & 124 & 230
		\\
		Rastrigin 5d & 159,564 & 180,539 & $-$ & 447 & $-$ & $-$
		\\
		Rastrigin 10d & 232,258 & 314,035 & $-$ & 1,004 & $-$ & $-$
	\end{tabular}
\end{table}

\end{document}